\documentclass{amsart}
\usepackage{amsfonts}
\usepackage{graphics}
\font\sr = cmr8

\newtheorem{theorem}{Theorem}
\theoremstyle{plain}

\newtheorem{claim}{Claim}

\newtheorem{corollary}{Corollary}

\newtheorem{lemma}{Lemma}

\newtheorem{proposition}{Proposition}
\newtheorem{remark}{Remark}

\numberwithin{equation}{section}

\begin{document}
\title[Acylindrical surfaces]{Acylindrical surfaces in 3-manifolds and knot complements}
\author{Mario Eudave-Mu\~noz}
\author{Max Neumann-Coto}
\address{Instituto de Matem\'aticas, UNAM\\ Ciudad Universitaria, 04510 M\'exico D.F., M\'exico. mario@matem.unam.mx max@matem.unam.mx}
\dedicatory{Dedicado a Fico en su 60 aniversario}

\begin{abstract}
We consider closed acylindrical surfaces in 3-manifolds and in knot and link complements, and show that the genus of these surfaces
is bounded linearly by the number of tetrahedra in a triangulation of the manifold and by the number of rational (or alternating)
tangles in a projection of a link (or knot). For each $g$ we find knots with tunnel number 2 and manifolds of Heegaard genus 3
containing acylindrical surfaces of genus $g$. Finally, we construct 3-bridge knots containing quasi-Fuchsian surfaces of unbounded
genus, and use them to find manifolds of Heegaard genus 2 and homology spheres of Heegaard genus 3 containing infinitely many
incompressible surfaces.
\end{abstract}

\maketitle

\section{Introduction}

A closed incompressible surface $F$ embedded in a 3-manifold $M$ is
called {\it acylindrical} if the manifold $M_{F}=M-intN(F)$, obtained by
cutting $M$
along $F$ contains no essential annuli (a properly embedded annulus in a 3-manifold is
essential if it is incompressible and not boundary parallel).
Acylindrical surfaces are interesting in connection
with geometry, as every totally geodesic surface in a hyperbolic 3-manifold
is acylindrical, and every acylindrical surface in a hyperbolic link
complement is quasi-Fuchsian. Moreover, if $F$ is an acylindrical surface in
a closed, irreducible and atoroidal 3-manifold $M$ then $M_{F}$ admits
a hyperbolic metric with totally geodesic boundary \cite{T}.

In \cite{H} Hass proved that for the finite volume hyperbolic 3-manifolds there
is a constant $C$, independent of the manifold, so that each acylindrical
surface in a manifold $M$ has genus at most $C\cdot vol(M)$. He
used this result to show that in any compact 3-manifold there is only a
finite number of acylindrical surfaces. It seems natural to ask if there are
similar bounds which hold for all 3-manifolds and depend not on
volume, but on some topological measures of complexity. Some candidates
could be the number of tetrahedra in a triangulation or the Heegaard genus
of the manifold, and in the case of knots and links, the crossing number, the bridge number or the tunnel
number. Such bounds must exist in the case of the number of tetrahedra in a
triangulation or the crossing number of a link, as there are only finitely
many manifolds and links for each number $n$. We find explicit bounds in these cases, and furthermore show that there is a linear bound in terms of the number of rational
tangles in a link projection or the number of alternating tangles in a prime knot projection.

The fact that 3-manifolds with Heegaard genus 2 and the complements of knots
with tunnel number 1 contain no separating acylindrical surfaces (\cite{MR}, \cite{BW})
could suggest that -at least for small Heegaard genus or tunnel number-
there could be bounds for the genus of such surfaces. We show
here that for each $g$, there are tunnel number 2 knots
which contain a closed acylindrical surface of genus $g$. By performing suitable Dehn
surgeries, we get closed manifolds of Heegaard genus $3$ which contain
closed acylindrical surfaces of genus $g$. These examples show that Heegaard genus
3 manifolds and tunnel number 2 knots are already quite complicated.

We also consider what happens when the acylindrical assumption is weakened to
require that there are no essential annuli running from the surface to
a boundary torus (in the case of hyperbolic knots and links this
means that the surface is quasi-Fuchsian).
We show that a knot that can be decomposed into two alternating tangles
cannot contain any quasi-Fuchsian surfaces in its complement.
On the other hand, we find hyperbolic 3-bridge knots whose complements
contain infinitely many quasi-Fuchsian surfaces. These knots have an
essential branched surface which carries quasi-Fuchsian
surfaces of arbitrarily high genus. These examples show that there are no bounds for
the genus of quasi-Fuchsian surfaces based on volume, crossing number or the
number of tetrahedra. Finally, by means of suitable Dehn fillings and double covers, we produce manifolds of Heegaard genus 2 and homology spheres of Heegaard genus 3 which contain infinitely many incompressible surfaces.
These examples are interesting, for it seems
that all known examples of hyperbolic manifolds with infinitely many surfaces
have noncyclic homology, and in the case of knots with infinitely many surfaces, it seems that the only known explicit examples are some satellite knots (see for example \cite{Ly}).
The examples are also interesting for the study of surfaces in the complement of 3-bridge knots, as they supplement results of Finkelstein and Moriah \cite{FM}, who showed that many 3-bridge knots contain an incompressible but meridionally compressible surface, and of Ichihara and Ozawa \cite{IO}, who proved that any closed surface in the complement of a 3-bridge knot is meridionally compressible or annular.

\section{Bounds for the genus of acylindrical surfaces}

\begin{proposition}
If a closed 3-manifold $M$ admits a (pseudo)triangulation with $n$
tetrahedra then the genus of a 2-sided closed acylindrical surface in $M$ is
at most $\frac{n+1}{2}$.
\end{proposition}

\begin{proof}
Let $T$ be a (pseudo)triangulation of $M$ with $n$ tetrahedra, and denote by
$T_{i}$ the $i$-skeleton of $T$.

Let $F$ be an incompressible surface in $M$ in normal position with respect
to the triangulation, so $F$ intersects the faces of the tetrahedra along
arcs and the interior of the tetrahedra along discs which are triangles or
squares. Assume further that $F$ has been isotoped to minimize the number of
intersections with $T_{1}$. Let $\overline{F}$ be the boundary of a regular
neighborhood $N$ of $F$. As $F$ is two-sided, $\overline{F}$ consists of two copies 
of $F$. By definition $F$ is acylindrical iff $M-intN$ contains no essential annuli.

The edges of $\overline{F}$ in each face of a tetrahedron split the face
into triangles, quadrangles, pentagons and/or hexagons, and each edge is
adjacent to a quadrangle (which lies in $N$). Call an edge {\it good} if
the other adjacent region (which lies in $M-intN$) is also a quadrangle.
Notice that if an embedded curve $c$ in $\overline{F}$ is made of good
edges, then the union of these adjacent quadrangles in $M-intN$ forms an
annulus $A$ that joins $c$ with another curve $c'$ in $\overline{F}$.

We claim that if $c$ is essential in $\overline{F}$ then the annulus $A$ is
essential. Otherwise $A$ would be isotopic to an annulus $A'$
bounded by $c$ and $c'$ in $\overline{F}$ (in particular, $c$ must
be 2-sided in $\overline{F}$). As $F$ is 2-sided in $M$, then $A'$
is parallel to an annulus $A''$ in $F$ and the isotopy from $%
A'$ to $A$ can be used to isotope $A''$ (pushing it
even further across $A$) to reduce the number of intersections of $F$ with $T_{1}$.

So if $\overline{F}$ is acylindrical, then the good edges of $\overline{F}$
carry no embedded essential curves, and so they carry no essential curves at
all. But as the edges of $\overline{F}$ split $\overline{F}$ into discs,
they must carry all of $H_{1}(\overline{F})$.

So there must be at least as many non-good edges in $\overline{F}$ as the
rank of $H_{1}(\overline{F})$. As the number of non-good edges in a face of
a tetrahedron is at most 6, the total number of non-good edges in
$\overline{F}$ is at most $12n$, so $12n\geq rank$ $H_{1}(\overline{F})=2\cdot genus$
$\overline{F}$, and so the genus of $F$ is at most $3n$.

In order to get the better estimate one needs to look more carefully at the
graph $Q$ formed by the edges of $\overline{F}$. Divide the
non-good edges in each tetrahedron in two classes: those lying in triangles
of $\overline{F}$ that cut off outermost corners of the tetrahedron will be
called {\it fair edges} and the others (which may lie in squares or
triangles) will be called {\it bad edges} (see Figure 1).

\bigskip
\centerline{\includegraphics{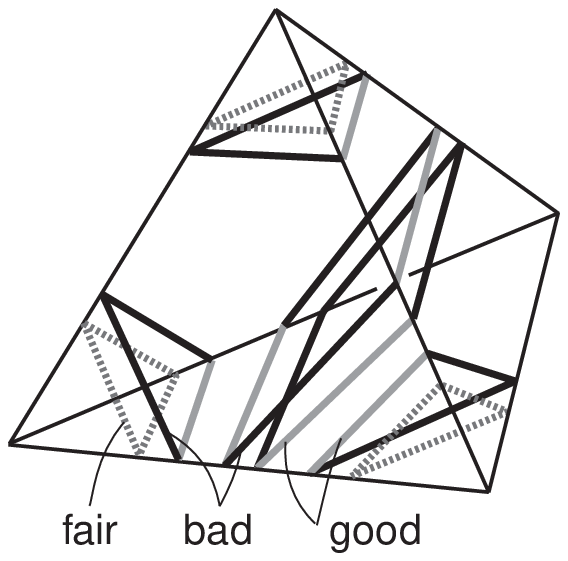}}
\centerline{\sr Figure 1}
\bigskip

Let $Q_g$ and $Q_f$ denote the subgraphs of $Q$ made of good edges and fair edges respectively.

Observe that $Q_g$ and $Q_f$ are disjoint, that is, have no vertices in common. As
the components of $Q_f$ lie in the links of the vertices of a
triangulation of the manifold $M$, all curves contained in $Q_f$ are contractible in
$M$, and so as $\overline{F}$ is incompressible then $Q_f$
contains only trivial curves of $\overline{F}$. All curves contained in $Q_g$ are also
trivial because $\overline{F}$ is acylindrical. So
as $Q_g\cup Q_f$ contains only trivial curves of $\overline{F}$, by
attaching to $Q_g\cup Q_f$
some of the complementary pieces of $\overline{F}$ we obtain a (possibly empty or
disconnected) simply connected subcomplex $\overline{F}_{S}$ of $\overline{F}$.

Now the Euler characteristic of $\overline{F}$ is $\chi (\overline{F})=\chi (
$ $\overline{F}_{S})+v-e+f$ where $v$, $e$, and $f$ count the vertices, edges
and faces of $\overline{F}$ that do not lie in $\overline{F}_{S}$. So $e$ counts
some bad edges -some others may lie in $\overline{F}_{S}$- and $f$ counts the 
triangles and squares adjacent to them. 
It can be shown directly that in each tetrahedron $\Delta $%
, every subcollection $f_{\Delta }$ of the set of squares and triangles of $%
\overline{F}$ $\cap \Delta $ with bad edges satisfies the inequality $\frac{1%
}{2}e_{\Delta }-f_{\Delta }\leq 2$. In particular, we may take $f_{\Delta }$
to be the set of squares and triangles with bad edges not contained in $%
\overline{F}_{S}$. As $\overline{F}$ has two components and each of them contains a
component of $\overline{F}_{S}$ or a vertex, it follows that

$\chi (\overline{F})\geq 2-e+f=2+\sum_{\Delta \in T_{3}}-\frac{1%
}{2}e_{\Delta }+f_{\Delta }\geq 2-2n$

and so $genus(F)=\frac{1}{4}rank H_1(\overline{F})=\frac{1}{4}(4-\chi(\overline{F}))
\leq \frac{1}{4}(2+2n)$.
\end{proof}

The genus of an acylindrical surface in a manifold is not bounded in terms
of its Heegaard genus, as we show in Section 3.  However, there is a bound depending on
the complexity of a Heegaard splitting. Let $M=H\cup H'$ be a Heegaard splitting of $M$
of genus $g$, and let $D_{1},D_{2},...,D_{g}$ and $D_{1}',D_{2}',...,D_{g}'$ be discs splitting $H$ and $H'$ into
3-balls $B$ and $B'$. The complexity of the Heegaard splitting with respect to these
discs is just the minimal intersection number between the boundaries of the discs.
The complexity of a Heegaard splitting is the minimum complexity among all such
systems of discs.

\begin{proposition}
If a closed 3-manifold $M$ admits an irreducible Heegaard splitting of
genus $g$ and complexity $n$ then the genus of a closed acylindrical
surface in $M$ is at most ($n-\frac{3}{2}g$).
\end{proposition}

\begin{proof}
Let $M=H\cup H'$ be a Heegaard splitting of $M$ of genus $g$ as above,
with $\cup D_{i}$ meeting $\cup D_{j}'
$ in $n$ points. Let $F$ be an acylindrical surface in $M$. As $F$ is
incompressible, we may assume that $F$ meets $H'$ along $g$ stacks of parallel discs in $N(D_{j}')$ (some stacks may be empty). We may also assume that $F$ meets $B$ along discs and that it meets each $D_{i}$ along stacks of parallel arcs connecting different
components of $\partial D_i \cap N(D'_j)$.

As before, consider the graph of intersection $Q$ of $\overline{F}$ with
$\partial H\cup _{i}D_{i}$.
Call an edge of $Q$ on $D_{i}$ {\it good} if it is an interior arc of a stack, otherwise call it {\it bad}. Call an edge of $Q$ in $\partial H$ {\it good} if it is part of the boundary of an interior disc of a stack. Otherwise (i.e., if it is part of the boundary of an outermost disc of a stack) call it {\it fair}. See Figure 2.

\bigskip
\centerline{\includegraphics{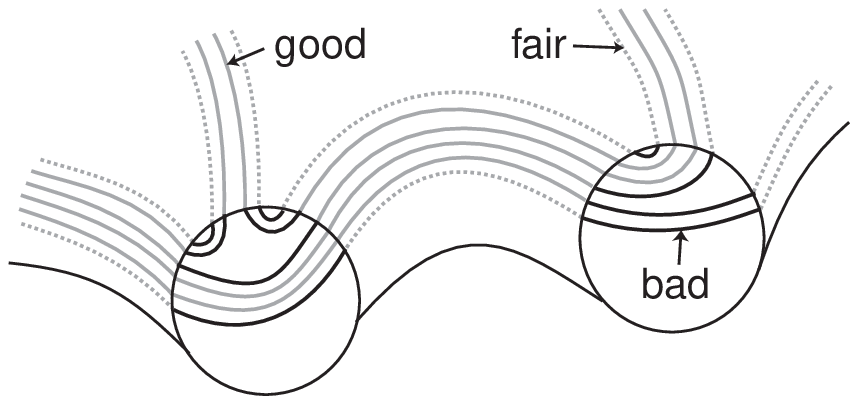}}
\centerline{\sr Figure 2}
\bigskip

Observe that the subgraphs $Q_g$ and $Q_f$ made of
good edges and fair edges do not meet. As the components
of $Q_f$ are contained in the boundaries of discs in $%
H'$ then $Q_f$ carries no essential curves of $\overline{F}$, and as $F$ is acylindrical we may assume as in the proof of 2.1 that $Q_g$ carries no essential curves either.

So as $Q_g \cup Q_f$ carries no essential curves and $Q$ splits $\overline{F}$
into discs, the rank of $H_{1}(\overline{F})$ is bounded above by the number of
bad edges. If $D_{i}$ meets the $D_{j}'$ in $n_{i}$ points then $n_{i}>1$
(because the Heegaard splitting is irreducible) and $D_{i}$ contains at most
$4n_{i}-6$ bad edges, so the rank of $H_{1}(\overline{F})$ is at most
$\sum_{D_{i}}\left( 4n_{i}-6\right)=4n-6g$
and so the genus of $F$ is at most $n-\frac{3}{2}g$.
\end{proof}

There are other ways of measuring the complexity of a Heegaard splitting,
for example, by means of the curve complex, as defined in \cite{He}.
Note however that no such bound for the genus of acylindrical surfaces
exists for this complexity, for in fact, all the examples constructed in Section 3 have a Heegaard splitting
of genus 3 which comes from a certain bridge presentation of a knot, and then by a similar proof to theorem 1.4 of
\cite{He}, the distance in the curvecomplex is $\leq 2$.

We now consider bounds for the genus of acylindrical surfaces in the exterior of knots and links in the 3-sphere.

\begin{proposition}
If $k$ is a knot or link with $n$ crossings then the genus of a closed
acylindrical surface in the exterior of $k$ is at most $\frac{3}{2}n-3$.
\end{proposition}

\begin{proof}
Draw $k$ on a projection sphere $S$, except for the crossings which lie on the
surface of $n$ small spheres $S_{1}$, $S_{2}$,...,$S_{n}$. Let $S_{0}$ be
the part of the projection sphere outside the $S_{i}$'s. Then $S_{0}\cup
_{i}S_{i}$ cuts $S^{3}$ into $n+2$ polyhedral balls $B^{-}$, $B^{+}$ and $%
B_{1}$, $B_{2}$, ..., $B_{n}$, with faces determined by the equators of the
bubbles and the arcs of $k$. If $F$ is an incompressible surface in the
exterior of $k$ then $F$ can be isotoped to meet $B^{+}$ and $B^{-}$ along
discs, meet each $B_{i}$ along parallel saddle-shaped discs, and meet their
faces along arcs. See Figure 3.

\bigskip
\centerline{\includegraphics{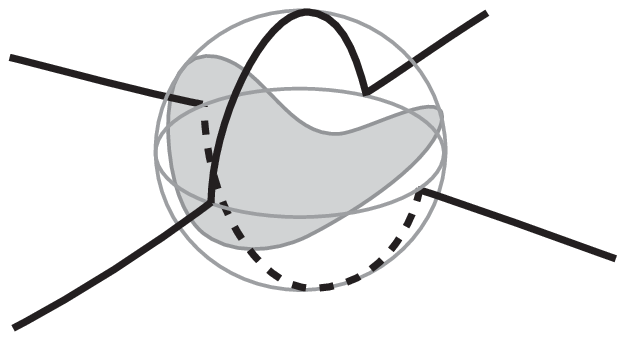}}
\centerline{\sr Figure 3}
\bigskip

Let $\overline{F}$ be the boundary of a regular neighborhood of $F$, and let
$Q$ be the graph of intersection of $\overline{F}$ with $%
S_{0}\cup _{i}S_{i}$. So $Q$ splits $\overline{F}$ into
discs. As before, consider the edges of $Q$ on each face of $S_{0}\cup _{i}S_{i}$,
call those that have parallel edges on both sides
{\it good}, those which are closest to arcs of $k$ and are parallel to
them {\it fair}, those lying on some $S_i$ and parallel to an arc of $\partial S_0$
which contain a point of $k$ are also fair,  and all the others are called
{\it bad} edges (so the faces of the $B_{i}$'s contain no bad edges). See Figure 4.

\bigskip
\centerline{\includegraphics{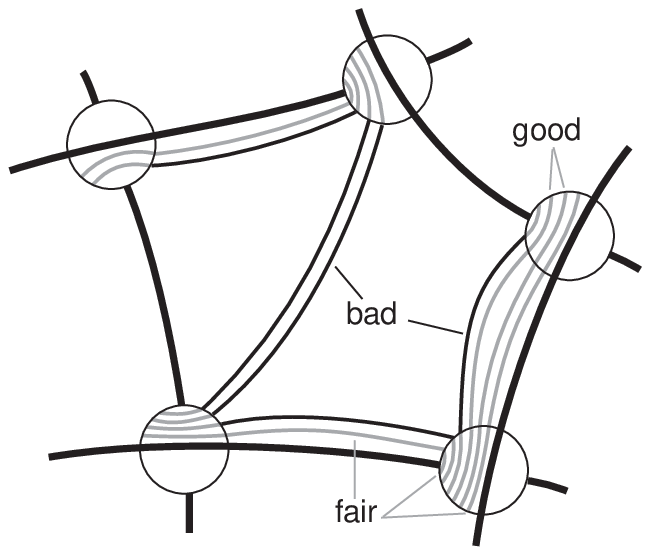}}
\centerline{\sr Figure 4}
\bigskip

Again the subgraphs $Q_g$ and $Q_f$ of $Q$
are disjoint, and if $\overline{F}$ is acylindrical then $Q_g$ carries no
essential curves of $\overline{F}$. On the
other hand, $Q_f$ can be regarded as lying on the boundary
tori of a regular neighborhood of the link $k$. But as $\overline{F}$ is
acylindrical, there can be no essential annuli running from $\overline{F}$
to $k$, so $Q_f$ contains no essential curves of $\overline{F}$. So, as the
graph $Q$ carries all of $H_{1}(\overline{F})$,
there must be at least as many bad edges as the rank of $H_{1}(\overline{F})$.
As there are at most $3i-6$ bad edges on each $i$-gon determined by the
projection of $k$ into $S$, the number of bad edges in $\overline{F}$ is at
most

\noindent $\sum_{i-\text{gons in }P}\left( 3i-6\right) =3(2($ arcs of $k$ in $S))-6\left(\text{regions determined by }k\text{ in }S\right) $

\noindent $=12n-6(2+n)=6n-12$

So the rank of $H_{1}(\overline{F})$ is at most $6n-12$ and the genus of $F$
is at most $\frac{6}{4}n-\frac{12}{4}$.
\end{proof}

After we proved proposition 3, we learned that Agol and D. Thurston,
following Lackenby \cite{L}, showed that the volume of a hyperbolic knot
of link is bounded above by $10 v_3(t(D)-1)$ where $v_3$ is the volume
of a hyperbolic ideal tetrahedra and
$t$ is the twist number of $k$ (the minimum number of
twists in a diagram of $k$, where a twist is a string of 2-gons or a
crossing in the diagram). Agol has also shown \cite{A} that if a hyperbolic
manifold
$M$ has an acylindrical surface of genus $g$,
then $Vol(M)\geq 4v_3(g-1)$. It follows
that the genus of an acylindrical surface in the exterior of a
hyperbolic link $k$ is at
most $\frac{5}{2}t$. These results suggested the following.

Recall that a tangle is a 3-ball $B$ together with two
properly embedded arcs. The tangle is rational if the arcs are isotopic
(rel $\partial$) to arcs in $\partial B$.
We will say that a knot or link $k$ in $S^3$ is decomposed into tangles if there is a
sphere $S$ and 3-balls $B_1$, $B_2$, ..., $B_n$ each intersecting
$S$ in a disc, so that $k\cap B_i$ is a tangle, and
the part of $k$ outside these balls is a
collection of arcs lying on $S_0=S-int(\cap B_i)$.

\begin{theorem}
If a link is decomposed into $n$ rational tangles, then the genus of a
closed acylindrical surface in its complement is at most $2n-4$.
\end{theorem}

\begin{proof}
Draw the projection of the link $k$ as the union of $n$ rational tangles
in the interior of $n$ disjoint spheres $S_{1}$, $S_{2}$,...,$S_{n}$ joined
by $2n$ disjoint arcs in the projection sphere. Let $S_{0}$ be the part of the
projection sphere outside these spheres. Then $S_{0}\cup _{i}S_{i}$ cuts $%
S^{3}$ into polyhedral balls $B_{+}$, $B_{-}$ and $B_{1}$,$B_{2}$,...,$B_{n}$
with faces determined by the equators of the spheres and the arcs of $k$ in $%
S_{0}$.

If $F$ is an incompressible surface in the exterior of $k$ we may isotope $F$ so that
intersects $B^{-}$ and $B^{+}$ along discs, and intersects $S$ and the
hemispheres of each $S_{i}$ along arcs. Moreover, as $k\cap B_{i}$ is a
rational tangle, we may isotope $F$ to intersect $B_{i}-k$ along parallel
discs that separate the strings of the tangle, and we may assume that their
boundaries meet each hemisphere of $S_{i}$ along 2 or 3 families of parallel
arcs -2 if the tangle is a crossing and 3 otherwise (a single family of
parallel arcs implies that the discs are vertical and the tangle has no
crossings of $k$).

Let $\overline{F}$ be the boundary of a regular neighborhood $N$ of $F$. The
intersection of $\overline{F}$\ with $S_{0}\cup _{i}S_{i}$ gives a cell
decomposition of $\overline{F}$ and cuts the faces of $S_{0}\cup _{i}S_{i}$
into quadrangles that lie in $N$ and other polygons that lie in $S^{3}-N$; as
before, let $Q$ be the graph of intersection. Call an
edge of $Q$ in a face of $S_{0}\cup _{i}S_{i}$ {\it good} if
the adjacent polygon in $S^{3}-N$ is a quadrangle with another edge on
$\overline{F}$ (so the two edges are parallel in that face). Otherwise, call
an edge in $Q$ {\it fair} if it is adjacent to a
quadrangle in $S_{0}$ with a side in $k\cap S_{0}$ that is adjacent to
another quadrangle in $S_{0}$ with a side in $Q$ (so both
edges of $Q$ are parallel to this arc of $k$) or if it is
adjacent to a polygon in a hemisphere of $S_{i}$ with exactly 2 sides in
$Q$ (so the other sides lie in the equator and are separated
by points of $k\cap S_{i}$). Call the other edges of $Q$
{\it bad}.  Note that edges lying on some $S_i$ and parallel to an arc of
$\partial S_0$
which contain a point of $k$ are bad. See Figure 5.

\bigskip
\centerline{\includegraphics{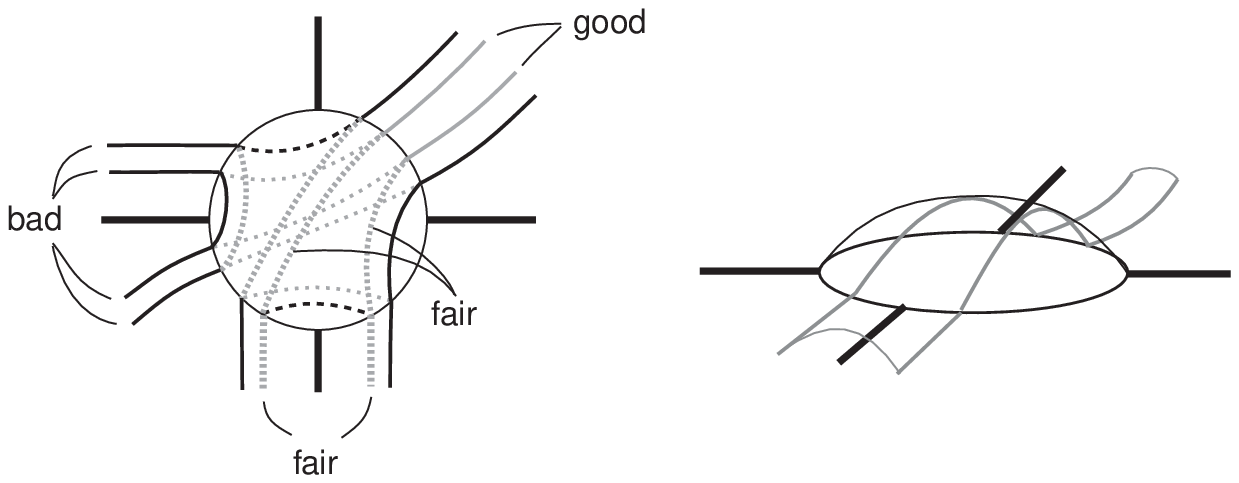}}
\centerline{\sr Figure 5}
\bigskip

One can use the fair edges as well as the good edges to construct annuli for
$\overline{F}$, by taking the quadrangles that lie between two fair edges in
$S_{0}\cup _{i}S_{i}$ (but that may intersect $k$) and pushing them outside
the corresponding $S_{i}$, or if they lie in $S_{0}$, to the side of $S_{0}$
that doesn't contain an edge of $Q$ connecting the two fair edges
(there can't be connecting edges on both sides because
the union of the four edges would be a meridian
of $k$, and so $F$ would be meridionally compressible). This creates
quadrangles in $S^{3}-k$\ connecting pairs of fair edges, and one can see
that the quadrangles corresponding to consecutive fair or good edges match
well.

As before, if a simple essential curve in $\overline{F}$ is made of good and
fair edges then the annulus formed by the union of the adjacent quadrangles
is essential or else $F$ could be isotoped to reduce its intersection with $%
S_{0}\cup _{i}S_{i}$. So, if $F$ is acylindrical, the subgraph of $Q$
consisting of the good and fair edges cannot contain any essential
curve of $\overline{F}$, so it is contained in a simply connected subcomplex
$\overline{F}_{S}$ of $\overline{F}$. Again, as $\overline{F}$ has two components and $Q$ divides them into discs, \break
$\chi (\overline{F})=\chi (\overline{F}_{S})+v-e+f \geq 2-e+f$ where
$v$, $e$, and $f$ count the vertices, edges
and discs in $\overline{F}-\overline{F}_{S}$, and so
$rank H_1(\overline{F})=(4-\chi(\overline{F}) \leq 2+e-f$.

As there are at most 4 bad edges and 8 fair edges on each $S_{i}$, all
contained in the 2 outermost discs of $\overline{F}\cap B_{i}$, the number of
bad edges minus the number of discs that contain them in $\cup S_{i}$ is at most $2n$.

There are at most $i-3$ families of parallel edges on each face of $S_{0}$
determined by $i>1$ arcs of $k$, not including the families of edges
parallel to the arcs of $k$, and they produce at most $2i-6$ bad edges on each face. If an arc of $k$ has parallel families
edges of $Q$ on both sides, then there are two bad edges in these families, for the edges closest to $k$ are fair. If an arc of $k$ has only edges of $Q$ on one side, then there are two bad edges in this family, since in this case the edge closest to $k$ is not fair.

So, if no face of $S_{0}$ is a
monogon the number of bad edges in $S_{0}$ is at most
$\sum_{edges of k} 2 + \sum_{i-\text{gons in
}P}\left( 2i-6\right)=4n + 8n-6(2+n)=6n-12$.

When $i=1$, the previous formula undercounts the number of bad edges in the
monogon as $-3$ instead of $0$ -there are no edges in the monogon as they
could be isotoped into $B_{i}$ to eliminate two intersection curves of $%
\overline{F}$ with $S_{i}$-. In this case there cannot be bad edges around
the endpoints of the monogon in $S_{i}$ and so the discs of intersection of $\overline{F}$ with $S_{i}$ are vertical and the tangle is trivial -unless $\overline{F}$ does not meet $S_{i}$ at all, so there is an overcount on the
number of bad edges in $\cup S_{i}$ by at least $2$ and also on the number
of bad edges in the face of $S_{0}$ adjacent to the monogon. So the
previous bound also holds when some faces of $S_{0}$ are monogons.

Finally observe that since $\overline{F}$ has 2 components and each of them must meet $B_{+}$ and $B_{-}$, there must be at least 2 discs of
$\overline{F}-\overline{F}_{S}$ inside each of these balls.

So, $genus(F)=\frac{1}{4}rank H_1(\overline{F})\leq \frac{1}{4}(2+e-f) \leq
\frac{1}{4}(2+2n+(6n-12)-4)=2n-\frac{7}{2}$
\end{proof}

Consider a tangle as above, i.e., it is determined by the intersection of a 3-ball $B$ with a link $k$, so that $B\cap S$ is a disc, where $S$ is a projection sphere, and $k\cap \partial B$ consists of 4 points lying on $S$. We say that the tangle is alternating if its arcs can be isotoped, keeping $\partial B$ fixed, to have an alternating projection on the sphere $S$. Note that each rational tangle is alternating.
The next result extends Theorem 1 to allow alternating tangles.

\begin{theorem}
If a prime knot is decomposed into alternating tangles, $n$ of them rational,
then the genus of a closed acylindrical surface in its complement is at most $2n-4$.
\end{theorem}

The proof is based on the following:

\begin{claim}
Let $k$ be a nonseparable link or a knot and $S$ a sphere that meets $k$ in 4 points.
Then each acylindrical surface $F$ in $S^3-k$ is isotopic to one that either
i) is disjoint from $S$ or
ii) intersects $S$ in one curve or
iii) meets one of the components of $S^3-k-S$ along parallel discs.
\end{claim}

\begin{proof}
The sphere $S$ separates $k$ into two tangles.
Isotope $F$ to minimize its intersection with the 4-punctured sphere $S-k$.
The intersection then contains no trivial curves, and as $F$ is meridionally
incompressible then it does not contain curves surrounding only one puncture, so
all the curves $c_1$,$c_2$,...$c_n$ in which $F$ intersects $S$ must be parallel in $S-k$.
As $F$ is acylindrical, if there is more than one $c_i$ then the annuli connecting two of
them in $S$ cannot be essential, so either one annulus is isotopic (rel $\partial $) to an
annulus in $F$ (and the isotopy can be used to remove two $c_i$'s) or all the $c_i$'s bound
discs of $F$. So at least one of them, say $c_1$, bounds a disc $D_1$ in $F$ that lies
completely on one side of $S$. But then, as all $c_i$'s are parallel to $\partial D_1$, one
can draw parallel discs $D_i$ in $S^3-k$ on that side of $S$ that meet $F$ at $c_i$ (and
nowhere else).
The union of the discs bounded by the $c_i$'s in $F$ and the $D_i$'s form spheres in $S^3-k$,
and if $k$ is a knot or a nonseparable link these spheres bound balls in $S^3-k$, so the
$D_i$'s must be isotopic to the discs in $F$, and the isotopy reduces the number of curves
unless the discs in $F$ were already on one side of $S$.
\end{proof}

\begin{claim}
If $k$ is a prime knot and $k\cap B_{i}$ is an alternating tangle, then every
acylindrical surface in the complement of $k$ can be isotoped to meet $%
B_{i}-k$ along parallel discs or be disjoint from it.
\end{claim}

\begin{proof}[Proof of theorem]
Assume for the moment that claim 2 is true, and isotope the surface $F$ to
meet only the $B_{i}$'s corresponding to separable tangles. To estimate the
genus of $F$ we would like to count the number of bad edges and discs of $%
\overline{F}$ that contain them by replacing each nonseparable tangle in the
diagram of $k$ by a trivial tangle to get a knot $k'$ and counting
the bad edges of $\overline{F}$ in its diagram.

\bigskip
\centerline{\includegraphics{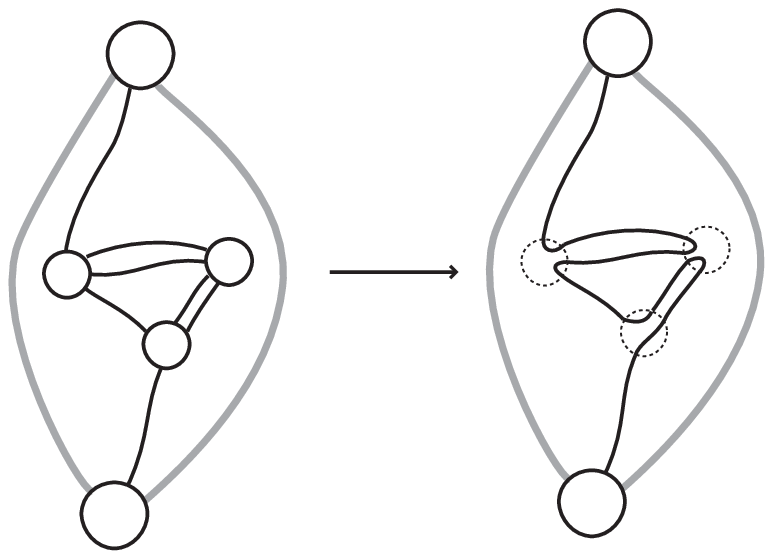}}
\centerline{\sr Figure 6}
\bigskip

Now some bad edges in the diagram of $k$ may become fair in the diagram of $k'$
as in Figure 6, but in this case we may regard them as originally
being ''almost fair'' -there is a quadrangle joining them that lies above or
below the nonseparable tangles that were between them in the diagram of $k$.
The quadrangles corresponding to almost fair edges match well with the other
quadrangles corresponding to good and fair pairs of edges, so they can be
used as well to construct annuli for $\overline{F}$. So the same bound for
the number of bad edges and discs -and therefore the same bound for the
genus of $F$- holds.
\end{proof}

The proof of claim 2 is based on the following extension of the Meridional Lemma of
Menasco \cite{M}.

\begin{lemma}
If a link $k$ intersects a ball $B$ in an alternating tangle, then every
meridionally incompressible surface in the complement of $k$ can be isotoped
to intersect $B$ along copies of a surface that separates the strings of the
tangle.
\end{lemma}

\begin{proof}
Draw $B$ as a round ball with $k\cap B$ lying in an equatorial disc except
at the crossings, that lie on the surface of small ``bubbles'' $B_{1}$, $%
B_{2}$,....as in Figure 3. Let $\partial B_{i+}$ and $\partial
B_{i-}$ be the hemispheres of $\partial B_{i}$, and let $D_{0}$ denote the
part of the equatorial disc outside the bubbles. Let $D_{+}=D_{0}\cup
_{i}\partial B_{i+}$ and $D_{-}=D_{0}\cup _{i}\partial B_{i-}$ , and let $%
B_{+}$ and $B_{-}$ be the parts of $B$ above and below $D_{+}$ and $D_{-}$.

If $F$ is a meridionally incompressible surface in the complement of $k$
then by isotoping $F$ to minimize its intersection with $\partial B\cup
D\cup _{i}\partial B_{i}$ we can assume that $F$ meets $\partial B$ along
parallel curves that separate 2 points of $\partial B\cap k$ from the other
2, that $F$ meets $D$ and each hemisphere of $\partial B$ and $\partial B_{i}
$ along arcs and to meet $B_{+}$ and $B_{-}$ along discs and each $B_{i}$
along parallel saddle-shaped discs. So $F$ intersects $D_{+}$ and $D_{-}$
along curves and arcs with endpoints in $\partial B$.

Following Menasco, one can show that the curves and arcs of intersection of $%
F$ with $D_{+}$ (and similarly with $D_{-}$) have the following properties:

1. As $F$ is incompressible, each curve (and each arc) crosses at least one
bubble.

2. As $F$ is meridionally incompressible, each curve (or arc) crosses each
bubble at most once.

3. As the diagram of $k\cap B$ is alternating, if a curve (or arc) crosses
two bubbles $B_{i}$ and $B_{j}$ in succession, then the 2 arcs $k\cap
\partial B_{i+}$ and $k\cap \partial B_{j+}$ lie on opposite sides of the
curve. See Figure 7.

\bigskip
\centerline{\includegraphics{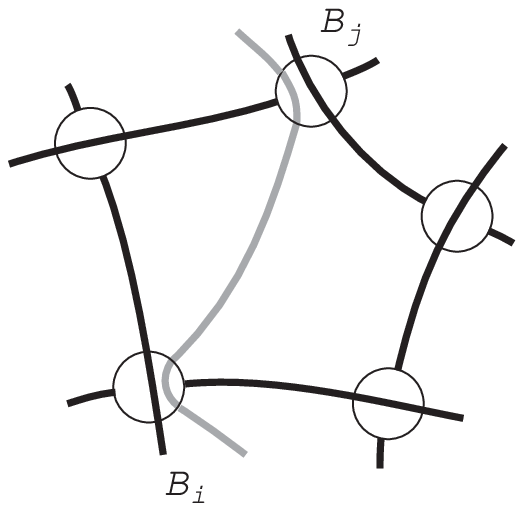}}
\centerline{\sr Figure 7}
\bigskip

So there can be no closed curves in $D_{+}$, because by properties 1 and 3
an innermost such curve would have to leave an arc of $k\cap \partial B_{i+}$
inside (so there would be another curve inside) unless the curve crossed the
same bubble twice, contradicting property 2.

Let $k^{1}$ and $k^{2}$, $k_{1}$ and $k_{2}$ be the 4 segments of $k \cap D_{0}$
that start on $\partial D_0$, and end in overcrossings or undercrossings of $k$
respectively. Note that $\partial D_0$ encounters them in the order $k^{1}$, $k_{1}$, $k^{2}$, $k_{2}$, for otherwise there is an arc on $D_{0}$ separating the strings of the tangle, but then the knot will be composite.
Properties 1, 2 and 3 for arcs imply that each outermost arc in $D_{+}$ goes
around $k_{1}$ or $k_{2}$ and so every
arc in $D_{+}$ must separate $k_{1}$ from $k_{2}$. See Figure 8a.

Now let $F_{0}$ be a surface consisting of one or more components of $F\cap B$. If $F_{0}$ does not separate the strings of the tangle then each path in $B$ joining the strings must meet $F_{0}$ in an even number of points, so $F_{0}$ intersects each bubble in an even number of discs, and so the number of curves and arcs cross $\partial B_{i+}$ on each side of $k\cap \partial B_{i+}$ is even. We claim that in these conditions $F\cap D_{+}$ consists of pairs of parallel arcs.

To show this, order the arcs according to its distance from $k_{1}$, and
assume that the first $2n$ are paired and let $a$ be the next one. Let $B_{i}$ and $B_{j}$ be two consecutive bubbles crossed by $a$, so the segments of $k\cap \partial B_{i+}$ and $k\cap \partial B_{j+}$ are on opposite sides of $a$ as in
Figure 8a. Since all the curves on one side of $a$ are paired and
each side of the bubbles is crossed by an even number of arcs, there must be
other arcs $a'$ and $a''$ crossing $B_{i+}$ and $%
B_{j+}$ next to $a$. If $a'$ and $a''$ are
different, then one of them cannot separate $k_{1}$ from $k_{2}$ (see Figure 8b).
If $a'=a''$ then either $a$
and $a'$ run parallel from $B_{i}$ to $B_{j}$ or else $a'$
crosses other bubbles between $B_{i}$ and $B_{j}$. If so, let $B_{l}$ be
the bubble crossed by $a'$ immediately after $B_{j}$. See Figure 8c.
Then $k\cap \partial B_{l+}$ lies between $a$ and $a'$, and so there must be another arc between $a$ and $a'$, and this
arc would have to cross $B_{i}$ or $B_{j}$ between $a$ and $a'$,
and this is impossible. Therefore $a'$ must run parallel to $a$ from the first
bubble to the last bubble crossed by $a$.
It remains to show that $a'$ runs parallel to $a$ from the first bubble to the
boundary of $D_{+}$ and from the last bubble to the boundary of $D_{+}$, i.e.,
that $a'$ does not meet other bubbles in its way to the
boundary and that the region between $a$ and $a'$ does not contain other bubbles.
As $k_1$ and $k_2$ lie outside the region between $a$ and $a'$,
this region does not contain any other arc $a''$. So $k^1$ and $k^2$ also lie outside this
region, because if $k^1$ were between $a$ and $a'$ the number of arcs between
$k_1$ and $k^1$ would be odd, so $F_{0}$ would separate these
strings of $k$. Now if there were any segments of $k \cap D_{+}$ in that region,
$k$ would have to enter and leave the region at 2 bubbles crossed by $a'$ on
its way to the boundary. But we know that for any two consecutive bubbles crossed by $a'$
the segments of $k$ in their upper hemispheres lie on opposite sides of $a'$,
so one of them is in the region between $a$ and $a'$ and so there must be
an arc in that region, a contradiction.

\bigskip
\centerline{\includegraphics{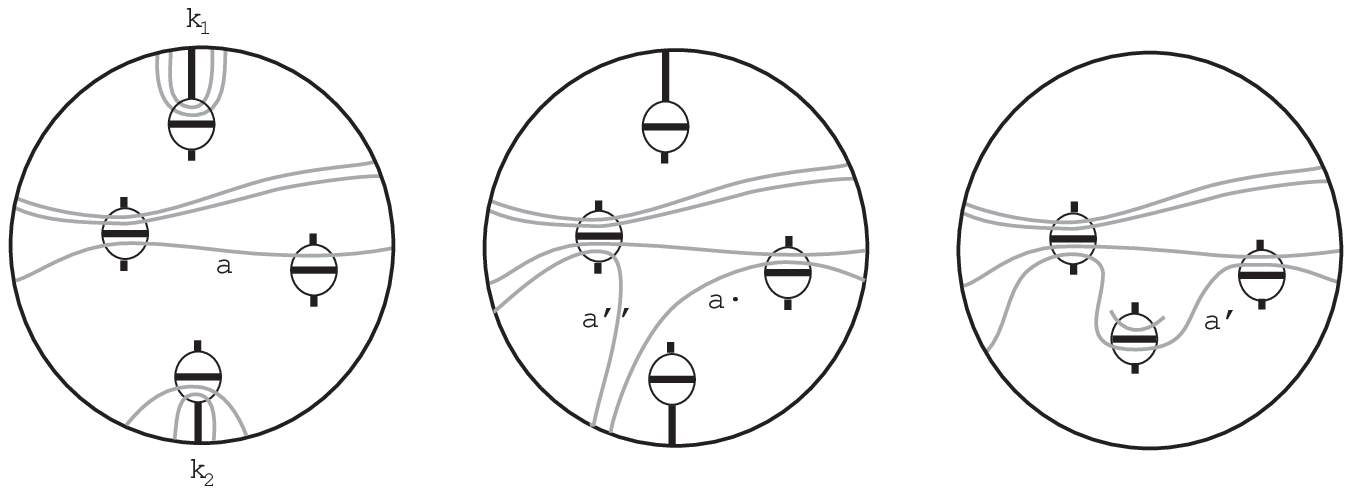}}
\centerline{\sr Figure 8}
\bigskip

Now observe that each pair of parallel arcs of $F_{0}$ in $D_{+}$ must be
adjacent to a pair of parallel arcs of $F_{0}$ in $\partial B_{+}$: an arc of
$F_{0}$ in $\partial B_{+}$ cannot go around the endpoints of $k^{1}$ or $
k^{2}$ because $F$ would be meridionally compressible  and something analogous holds
for the arcs of $F_{0}$ in $D_{-}$. So the intersection of $F_{0}$ with
$\partial B_{+}$, $\partial B_{-}$ and with each $\partial B_{i}$ consists of
pairs of parallel curves, and as $F_{0}$ is assembled by attaching discs to
these parallel curves, $F_{0}$ must consist of pairs of parallel surfaces.

Finally, as $F$ is meridionally incompressible, the intersection of $F$ with
the 4-punctured sphere $\partial B-k$ consists of curves surrounding 2
punctures, and if there is more than one curve these are parallel. So if $%
F\cap B$ has several components, and $F_{0}$ consists of any two of them,
then any path in $B$ joining the strings of the tangle must intersect $F_{0}$
in an even number of points, and this is all that we needed before to show
that $F_{0}$ consists of parallel surfaces.
\end{proof}

\begin{proof}[Proof of claim 2] Isotope $F$ to minimize its intersection with $\partial B_{i}$.
By the previous lemma if $F\cap B_{i}$ is
not empty then it consists of parallel copies of a surface $F_{0}$ that separates
the strings of the tangle. As $k$ is a knot $F$ cannot separate the strings
of $k\cap B$, so there must be an even number of copies of $F_{0}$.
Now by the previous claim either $F\cap B_{i}$ or $F\cap S^3-B_{i}$ consists of discs,
and in the second case $F$ would be the union of the components of $%
F\cap B_{i}$ with discs, and since there are at least two such components $F$ would not be connected.
\end{proof}

In \cite{ABB} Adams {\it et al.} extended the Meridional Lemma of Menasco to
almost alternating knots, i.e. knots that can be obtained by
changing one crossing of an alternating knot. The following corollary
extends it to knots that can be obtained from an alternating one by
mirroring any (2-string) tangle.

\begin{corollary}
If a knot $k$ can be decomposed into 2 alternating tangles, then $k$ admits
no meridionally incompressible surfaces in its complement.
\end{corollary}

\begin{proof}
By lemma 1, a meridionally incompressible surface $F$ in the complement of $k$
can be isotoped to meet each of the balls $B_{1}$ and $B_{2}$ that determine the tangles
along an even number of parallel copies of a surface $F_{i}$ that separates
the strings of the tangle.

So $F\cap B_{i}$ is the boundary of a regular neighborhood $N_{i}$ of one or
more copies of $F_{i}$, and $N_{i}$ is determined by painting the components
of $B_{i}-F$ in a chessboard fashion and choosing those whose color is
different from that of the regions that contain the strings of the tangle.
So $N_{1}$ and $N_{2}$ match on $\partial B_{1}=\partial B_{2}$ to form the
regular neighborhood of a single surface in $S^{3}$, and $F$ is its
boundary, so $F$ cannot be connected.
\end{proof}

\begin{corollary}
The total genus of a disjoint family of closed, embedded, totally geodesic surfaces in a
hyperbolic 3-manifold or link complement is bounded above by:

\begin{itemize}
\item $\frac{3}{2}t$ where $t$ is the number of tetrahedra in a triangulation.

\item $n-\frac{3}{2}g$ for manifolds of Heegaard genus $g$ and complexity $n$.

\item $\frac{3}{2}c-3$ for a link with $c$ crossings.

\item $\frac{5}{2}r-3$ for a link that admits a projection made of $r$ rational tangles.

\item $\frac{5}{2}r-3$ for a prime knot decomposed into alternating tangles, $r$ of them rational.

\end{itemize}
\end{corollary}

\begin{proof}
If $M$ is a hyperbolic 3-manifold and $F_1$,$F_2$,...,$F_k$
are disjoint totally geodesic surfaces in $M$, then each $F_i$ is
acylindrical and there are no essential annuli in $M$ connecting two $F_i$'s.
For, the preimages of the $F_i$'s in the universal covering of $M$
are disjoint totally geodesic planes in $\mathbb{H}^3$, and each preimage of an
essential annulus is an infinite strip of bounded height connecting
two lines in different planes.
These lines lie at a bounded distance from geodesic lines representing
the preimages of the boundaries of the annulus, so they determine 2 different
points at infinity where the two planes meet, but two disjoint totally
geodesic planes in $H^3$ can only meet at 1 point.

So we may consider the family $F_1$,$F_2$,...,$F_k$ as a single disconnected
acylindrical surface.
The arguments above show the existence of essential annuli for a surface $F$
if the rank of $H_1(F)$ is higher than the number of bad edges, independently of the
number of components of $F$. The bounds arise from a count of the number of
bad edges in each case.
\end{proof}

\section{Acylindrical surfaces in tunnel number two complements}

Let $S$ be a closed surface of genus $g$ standarly embedded in $S^3$,
that is, it bounds a handlebody on each of its sides. A knot $K$ has a
$(b,g)$-presentation if can be isotoped to intersect $S$
transversely in $2b$ points that divide $K$ into $2b$ arcs, so that the
$b$ arcs in each side can be isotoped, keeping the endpoints fixed, to
disjoint arcs on $S$. We say that a knot $K$ is a $(b,g)$-knot if it
has a $(b,g)$-presentation. Consider a product neighborhood $S\times I$
of $S$. To say that a knot $K$ has a $(b,g)$-presentation is equivalent to
say that $K$ can be isotoped to lie in $S\times I$, so that
$K\cap (S\times \{ 0 \})$ and $K\cap (S\times \{ 1 \})$ consist each of
$b$ arcs (or $b$ tangent points), and the rest of the knot consist of
$2b$ straight arcs in $S\times I$, that is, arcs which intersect each
leave $S\times \{t\}$ in the product exactly in one point. It is not
difficult to see that if $K$ is a $(b,g)$-knot, then the tunnel number
of $K$, denoted $tn(K)$, satisfies $tn (K)\leq b+g-1$.
In this section we construct $(2,1)$-knots, which are in fact tunnel number 2 knots,
which contain an acylindrical surface of genus $g$.

Let
$T$ be a standard torus in
$S^3$, and let
$I=[0,1]$.  Consider $T\times I\subset S^3$.
$T\times \{ 0\}$ bounds a solid torus $R_0$, and
$T\times \{1\}$ bounds a
solid torus $R_1$, such that $S^3= R_0\cup (T\times I)\cup R_1$.
Choose $n+1$ distinct points on $I$, $e_0=0,\ e_1,\dots, e_n=1$, so that
$e_i < e_{i+1}$, for all $0\leq i\leq n-1$. Consider the tori
$T\times \{e_i\}$.
By a vertical arc in a product $T\times [a,b]$ we mean an
embedded arc which intersects every torus $T\times \{x\}$  in
the product in at most one point.

Let $\gamma_i$ be a simple closed essential curve embedded in the
product $T\times [e_{i-1},e_i]$, for $i=1,\dots,n$, so that it has only
one local maximum and one local  minimum with respect to the projection
to $[e_{i-1},e_i]$.
Let $\alpha_i$, for $i=1,\dots,n-1$, be a vertical arc in
$T\times [0,1]$, joining the maximum point of $\gamma_i$ with the
minimum of  $\gamma_{i+1}$. Let $\Gamma$ be the 1-complex consisting of
the union of all the curves $\gamma_i$ and the arcs $\alpha_j$. So
$\Gamma$ is a trivalent graph embedded in $S^3$. Let $R_0'=R_0\cup
(T\times[e_0,e_1])$ and $R_1'=R_1\cup (T\times[e_{n-1},e_n])$.

Suppose each curve $\gamma_i$ satisfies the following:

\begin{enumerate}
\item $\gamma_i$ is not in a 3-ball contained in
$T\times [e_{i-1},e_i]$, or in $R_0'$ or $R_1'$, that is, it is not a
trivial knot in that region.

\item  $\gamma_i$ is not isotopic in $T\times [e_{i-1},e_i]$, or in
$R_0'$ or $R_1'$, to a knot lying on the torus
$T\times \{e_i\}$.

\item $\gamma_i$ is not a cable of a knot lying in
$T\times [e_{i-1},e_i]$ or in
$R_0'$ or $R_1'$ (it can be proved that this is
equivalent to say that $\gamma_i$ is not isotopic to a cable of a knot
lying on the torus $T\times \{e_i\}$.)

\item
There is no  annulus $B$ in $T\times \{e_0\}$ so that
$B \times [0,1]$ contains $\Gamma$. If that happens then each
curve $\gamma_i$
would be contained in a product $B \times [e_{i-1},e_i]$.

\item There is no  M\"obius band in $R_0'$
($R_1'$) disjoint from $\gamma_1$ ($\gamma_n$).

\end{enumerate}

It is not difficult to see that there exist plenty of knots satisfying the
conditions required for the curves $\gamma_i$, say by taking each
$\gamma_i$ to be a $(1,1)$-knot which is not a torus knot nor a satellite
knot. For example, each
$\gamma_i$ could be a copy of the figure eight knot, as shown in Figure
9(a) in the case of $\gamma_1$, Figure 9(b) for
$\gamma_2,\dots,\gamma_{n-1}$, and Figure 9(c) for $\gamma_n$. In the
figures the knot is divided in two arcs; the thin arc contains the
minimum point of the knot, and the bold arc contains the maximum. When
assembled we get the graph $\Gamma$, shown for
$n=2$ in Figure 10.

\bigskip
\centerline{\includegraphics{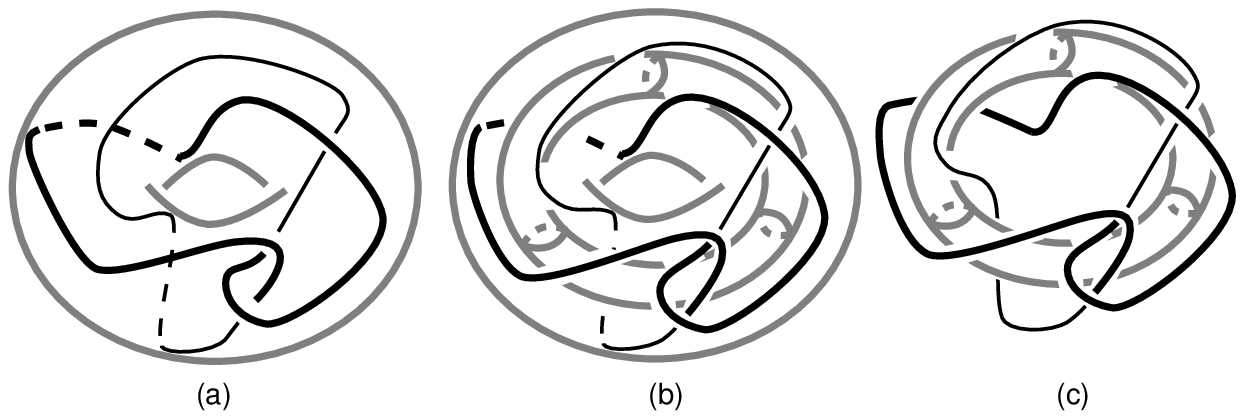}}
\centerline{\sr Figure 9}
\bigskip

\bigskip
\centerline{\includegraphics{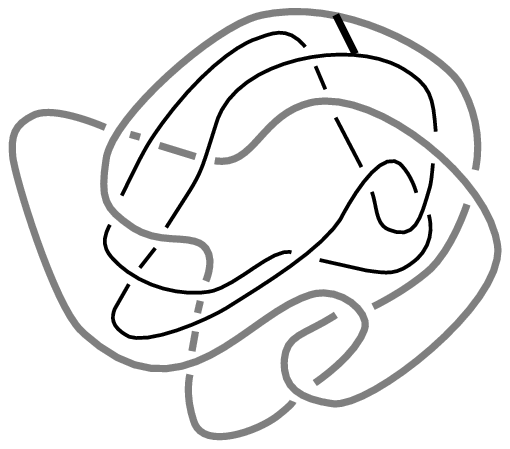}}
\centerline{\sr Figure 10}
\bigskip

Let $N(\Gamma)$ be a regular neighborhood of
$\Gamma$. This is a genus $n$ handlebody. We can assume that $N(\Gamma)$
is the union of $n$ solid tori $N(\gamma_i)$, joined
by $(n-1)$ 1-handles $N(\alpha_j)$.

\begin{theorem}
Let $\Gamma$ be a graph as above. Then
$S=\partial N(\Gamma)$ is incompressible and acylindrical in
$S^3-intN(\Gamma)$.
Furthermore, $M-intN(\Gamma)$ is atoroidal.
\end{theorem}

\begin{proof}
Consider the tori $T\times \{e_i\}$, $1\leq i \leq n-1$. These
tori divide $S^3$ into $n+1$ regions, where $n-1$ of them are product
regions and two of them are solid tori, namely
$R_0'$ and $R_1'$.
The torus
$T\times \{e_i\}$ intersects $\Gamma$ in one point, that is, a
middle point of $\alpha_i$, so $T\times \{e_i\} \cap N(\Gamma)$ consists
of a disc. Let $T_i = T\times \{e_i\} - int N(\Gamma)$, for
$1\leq i \leq n-1$, this is a once punctured torus.

Suppose $D$ is a compression disc for $S$, and suppose it intersects
transversely the tori $T_i$. Let $\beta$ be a simple closed curve of
intersection between $D$ and the collection of tori, which is innermost
in $D$. So $\beta$ bounds a disc $D' \subset D$, which is contained in a
product  $T\times [e_{i-1},e_i]$, or in the solid torus  $R_0'$
or in $R_1'$. If $\beta$
is trivial on $T_i$, then by cutting $D$ with an innermost disc lying in the disc bounded by $\beta$
on $T_i$, we get a compression disc with fewer intersections with the
$T_i's$. If $\beta$ is essential on $T_i$, then it would be parallel to
$\partial T_i$, or it would be a meridian of $T_1$ or a longitude of
$T_{n-1}$, but then in any case, one of the curves
$\gamma_1$ or $\gamma_n$ will be contained in a 3-ball, which is a
contradiction.

So suppose $D$ intersects the $T_i's$ only in arcs. Let $\beta$ such
an arc which is outermost on $D$, then it cobounds with an arc
$\delta \subset \partial D$ a disc $D'$. We can assume that $\beta$ is
an arc properly embedded in some $T_i$; if $\beta$ is parallel to an arc
on $\partial T_i$, then by cutting $D$ with an outermost such arc lying
on $T_i$ we get another compression disc with fewer intersections with
the $T_i's$, so assume that $\beta$ is an essential arc on $T_i$.
After isotoping $D$ if necessary, we can assume that the arc $\delta$ can
be decomposed as
$\delta = \delta_1 \cup \delta_2 \cup \delta_3$, where
$\delta_1, \delta_3$ lie on $\partial N(\alpha_i)$ and $\delta_2$
lie on $\partial N(\gamma_i)$ (if $\delta$ were contained in
$\partial N(\alpha_i)$, then by isotoping $D$ we would get a compression
disc intersecting $T_i$ in a simple closed curve). Let
$E$ be a disc contained in $N(\alpha_i)$ so that
$\partial E = \delta_1 \cup \delta_4 \cup \delta_3 \cup \delta_5$, where
$\delta_4$ lies on $T_i$ and $\delta_5$ lies on $\partial N(\alpha_i)$.
So $D' \cup E$ is an annulus, where one boundary component, i.e.,
$\beta \cup \delta_4$ lies on $T\times \{e_i\}$, and the other,
$\delta_2 \cup \delta_5$, lies on $\partial N(\gamma_i)$. If $\delta_2
\cup \delta_5$ is a meridian of $\gamma_i$, then necessarily $D\cup E$ is
contained  in $R_0'$ (or in $R_1'$)
and $\beta \cup \delta_4$ is a meridian of that solid torus. Then
$\gamma_1$ (or $\gamma_n$) intersects a meridian disc of $R_0'$ ($R_1'$) in one
point, which implies that it is parallel to a knot lying on the torus
$T\times \{e_0\}$ ($T\times \{e_1 \}$), which is a contradiction. If
$\delta_2 \cup \delta_5$ is a longitudinal curve of $\gamma_i$, then
this implies that $\gamma_i$ is parallel to a curve on
$T\times \{e_i\}$, a contradiction. If $\delta_2 \cup \delta_5$ goes
more than once longitudinally on $\gamma_i$, this would only be possible
for the curves $\gamma_1$ or $\gamma_n$, but then one of these curves
would be a core of the solid torus $R_0'$ or $R_1'$, which is not
possible.  This completes the proof that $S$ is
incompressible in $S^3-int N(\Gamma)$.

Suppose now that there is an essential annulus $A$ in
$S^3-int N(\Gamma)$. Look at the intersection between $A$ and the
punctured tori $T_i$. Simple closed curves of intersection which are
trivial on $A$, and arcs on $A$ which are parallel to a component of
$\partial A$ are eliminated as above. So the intersection consists of a
collection of essential arcs on $A$, or a collection of essential simple
closed curves on $A$.

Suppose first that there are essential arcs of intersection. Let
$E \subset A$ be a square determined by the arcs of intersection. So
$\partial E = \epsilon_1 \cup \delta_1 \cup \epsilon_2 \cup \delta_2$,
where $\epsilon_1$, $\epsilon_2$ are contained in different components
of $\partial A$ and $\delta_1$, $\delta_2$ are arcs of intersection of
$A$ with the $T_i's$. Take the square at highest level. So
$\delta_1$, $\delta_2$ lie on the same level $T_i$, and possibly
$T_i=T_{n-1}$. So we can assume that $\epsilon_1$, $\epsilon_2$ lie on
$\partial N(\alpha_i \cup \gamma_{i+1})$.

\medskip
{\sl Case 1: The arcs $\delta_1$, $\delta_2$ are parallel on $T_i$,
that is, they cobound a disc $F$ in $T_i$.}

There are two subcases, depending of the orientation of the arcs
$\delta_1,\delta_2$. Give an orientation to $\partial E$. Suppose first
that the arcs $\delta_1,\delta_2$ have the same orientation on $T_i$ (note
that the interior of $F$ may intersect the annulus $A$, but it is
irrelevant in this case). Then $E\cup F$
is a M\"obius band, and by pushing it off $T_i$ we get a M\"obius band
contained in the product $T\times [e_i,e_{i+1}]$ or in $R_1'$, with its
boundary lying on $N(\gamma_i)$. This implies
that either $\gamma_i$ is a trivial knot or that it is a 2-cable of some
knot, which is a contradiction.

Suppose the arcs $\delta_1,\delta_2$ have opposite orientations in $T_i$. If the interior of the disc $F$ intersects $A$, then take another square in $A$, which determines a disc $F'\subset F$ with interior disjoint from $A$.
We can form two annuli, $E\cup F$ and $(A-E)\cup F$.
We will show that at least one of them is an essential annulus.
Note that a core of $A$ is homotopic to the
product  of a core of $E\cup F$ and a core of $(A-E)\cup F$. So if these
two curves are homotopically trivial, so is  the core of $A$. So
assume one of them  is incompressible, say $(A-E)\cup F$. If it is
$\partial$-compressible then it is $\partial$-parallel, because $S$
is incompressible. Then there is a $\partial$-compression disc for this
annulus intersecting it on $(A-E)$, but this implies that the original
annulus $A$ is also $\partial$-compressible, a contradiction.
So we get a new essential annulus with fewer
intersection with the $T_i's$.

\medskip
{\sl Case 2: The arcs $\delta_1$, $\delta_2$ are not parallel on
$T_i$,
and the arcs $\epsilon_1$, $\epsilon_2$ are parallel on
$\partial N(\Gamma)$.}

The arcs $\epsilon_1$, $\epsilon_2$ must have the same
orientation on $N(\alpha_i\cup \gamma_i)$, see Figure 11(a). They cobound a
disc $F$ on
$\partial N(\alpha_i\cup \gamma_i)$ with
$\partial F = \epsilon_1\cup \eta_1\cup \epsilon_2\cup\eta_2$, where
$\gamma_1,\gamma_2 \subset \partial N(\alpha_i)\cap T_i$. (Note that the
disc $F$ may intersect the arc $\alpha_{i+1}$, or its interior may intersect $A$, but this is irrelevant in
this argument). It follows that
$E\cup F$ is a M\"obius band whose boundary lies on $T_i$. This is
impossible if the band lies in a product region. If it lies in $R_1'$,
then note that the band is disjoint from the curve
$\gamma_n$, but this is not possible, by hypothesis.

\bigskip
\centerline{\includegraphics{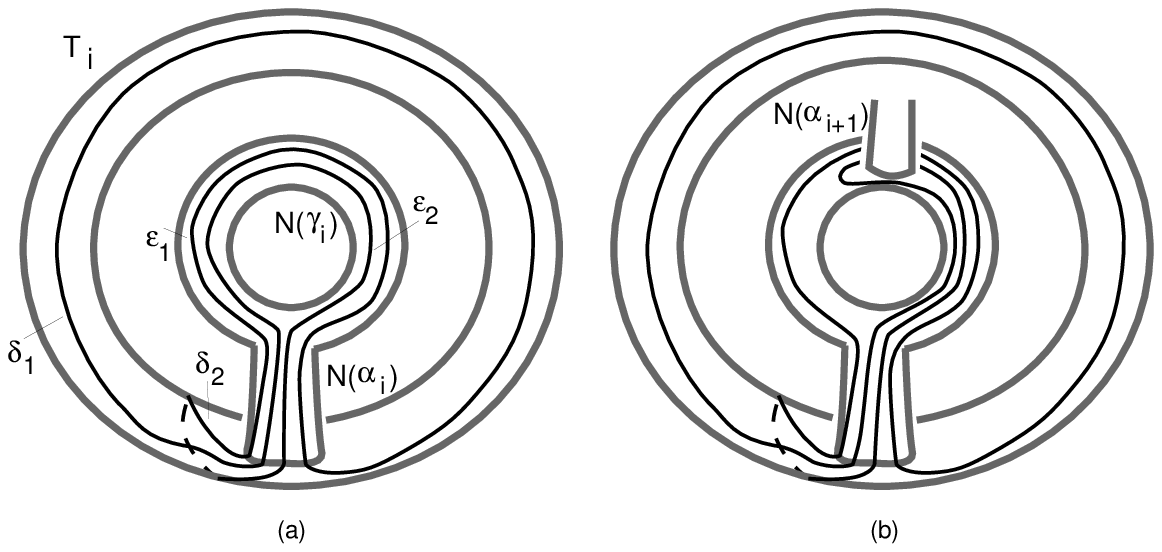}}
\centerline{\sr Figure 11}
\bigskip

{\sl Case 3: The arcs $\delta_1$, $\delta_2$ are not parallel on
$T_i$, and the arcs $\epsilon_1$, $\epsilon_2$ are not parallel on
$\partial N(\Gamma)$.}

Note that this case is only possible in a product region, see Figure 11(b).
Forget about the arc
$\alpha_{i+1}$, that is, consider the  square $E$ in the complement of
$N(\alpha_i\cup
\gamma_i)$. Then, it is not difficult to see that one of the arcs, say
$\epsilon_2$ can be slid toward $T_i$. Then there is a disc, whose
boundary consists of two arcs, one lying on $T_i$ and one on
$N(\gamma_i)$. By gluing to this disc a disc contained in $N(\alpha_i)$,
an annulus between
$\gamma_i$ and $T\times e_i$ is constructed. The only possibility in this
case is that the annulus goes once longitudinally on $N(\gamma_i)$, i.e.,
the curve $\gamma_i$ is parallel to the torus $T\times e_i$, which is a
contradiction.

This completes the proof in the case the annulus $A$ is divided
in squares.

Suppose now that the intersection of the annulus $A$ with the tori $T_i's$
consists of simple closed curves which are essential on $A$. Take an
outermost curve, say $\alpha$. Then $\alpha$ and a component of
$\partial A$ cobound an annulus, and the component of $\partial A$ must
lie on some $\gamma_i$. This again implies that $\gamma_i$ is parallel to
$T_i$ or that $\gamma_1$ or $\gamma_n$ are the core of the solid torus
$R_0'$ or $R_1'$, a contradiction.

It remains to prove that $S^3 -int N(\Gamma)$ is atoroidal. Suppose $Q$
is an essential torus, then we can assume that it intersects the tori $T_i$ in a
collection of simple closed curves which are essential on $Q$, and divide
$Q$ in a collection of annuli. Take one of this annuli, say $A$, at
highest level. If $A$ is in a product region then it must be parallel to
some $T_i$, and then by an isotopy we can remove two curves of
intersection. So $A$ lies on $R_1'$. As it is an annulus in a solid torus,
it must be parallel to the boundary. If $\gamma_n$ is not in this
parallelism region, then an isotopy removes the intersection. If
$\gamma_n$ is the parallelism region, then take the annulus next to $A$.
It must be an annulus between $T_{n-1}$ and $T_{n-2}$. Continuing in
this way, the only possibility is that the whole graph $\Gamma$ lies
inside a solid torus bounded by $Q$, but this is isotopic to a
solid torus of the form $B\times I$, where $B$ is an annulus
in $T\times \{ e_n\}$. This contradicts the choice of $\Gamma$.
\end{proof}

Put now a knot $K$ inside $N(\Gamma)$ in such a way that
$K \cap N(\alpha_i)$, for $2\leq i \leq n-1$, consists of four vertical
arcs with a pattern like in Figure 12(c), and
$k\cap N(\gamma_i)$ consists of 4 vertical arcs, going from $N(\alpha_i)$
to $N(\alpha_{i+1})$, as in Figure 12(b). Also, $K\cap N(\gamma_1)$ consists
of two arcs, each having a single minimum, and $K\cap N(\gamma_n)$ consists
of two arcs, each having a single maximum, as in the pattern shown in
Figure 12(a). For $n=3$, a knot $K$ inside $N(\Gamma)$ looks like
in Figure 13, where the twist is added to get a knot.

\bigskip
\centerline{\includegraphics{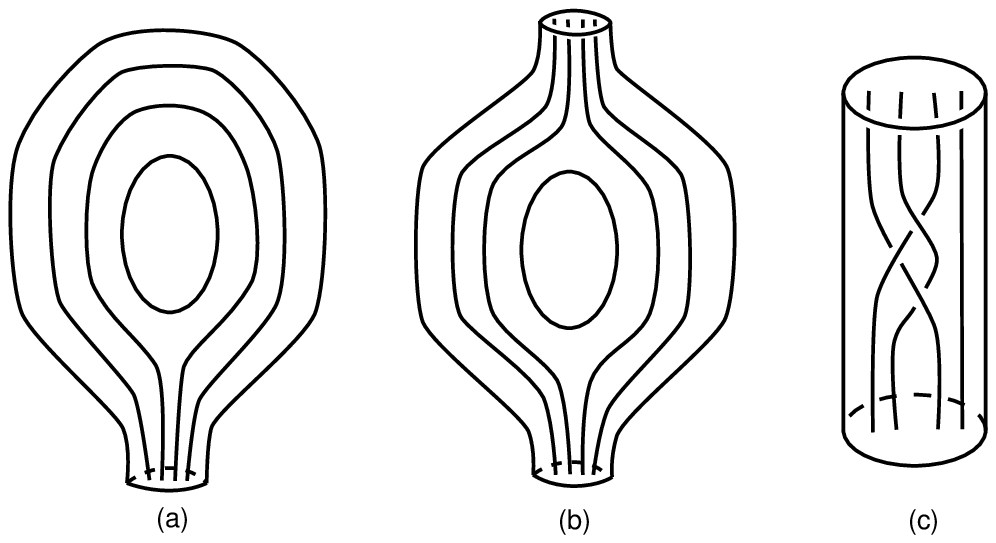}}
\centerline{\sr Figure 12}
\bigskip

\bigskip
\centerline{\includegraphics{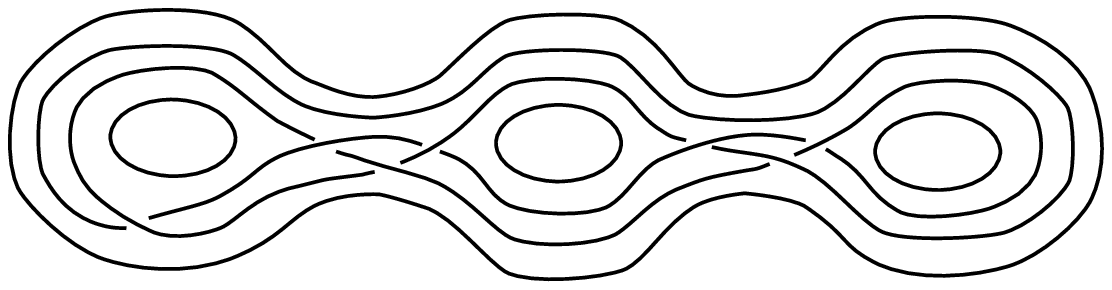}}
\centerline{\sr Figure 13}
\bigskip

\begin{lemma}
$S=\partial N(\Gamma)$ is acylindrical in
$N(\Gamma)-K$.
Furthermore, $N(\Gamma)-K$ is atoroidal.
\end{lemma}

\begin{proof} 
The proof is also an innermost disc/outermost arc
argument. It is practically the same as in Lemma 2.3 of \cite{AR}.
\end{proof}

\begin{theorem}
Let $K$ and $S$ as constructed above.
$K$ is a hyperbolic $(2,1)$-knot, tunnel number 2 knot, and $S$ is an acylindrical surface of genus $g$ in the complement of $K$.
\end{theorem}

\begin{proof}
Note that by construction $K$ is a $(2,1)$-knot, for it
lies in $T\times I$, and it has in there exactly two maxima and two
minima  with respect to the projection to the factor $I$.
It follows from Theorem 3 and lemma 2 that $S$ is an acylindrical surface.
$K$ is a hyperbolic knot because the complement of the surface is atoroidal and acylindrical.
Finally note that the knot $K$ has tunnel number 2; it cannot have tunnel
number one, for it contains an acylindrical separating surface \cite{MR}.
\end{proof}

\begin{corollary} Given any integer $g\geq 2$, there exist
infinitely  many hyperbolic 3-manifolds of Heegaard genus 3 which
contain an acylindrical surface of genus $g$.
\end{corollary}

\begin{proof}
For each $g$ choose a knot $K$ as above. Do Dehn
surgery on $K$ with slope $\lambda$, such that $\Delta(\mu,\lambda)\geq
3$, where $\mu$ is a meridian of $K$. It follows that $S$ remains
incompressible \cite{W}, acylindrical \cite{GW},
and that $M(\alpha)$ is irreducible and atoroidal \cite{GL1} \cite{GL2}.
Then by Thurston Geometrization Theorem, $M(\alpha)$ is hyperbolic, for it
is Haken and atoroidal. $K$ has tunnel number two, which implies that
$M(\alpha)$ has Heegaard genus at most 3, but it cannot have Heegaard
genus 2, for it contains a separating acylindrical surface \cite{MR}.
\end{proof}

\section {Quasi-Fuchsian surfaces of arbitrarily high genus}

Let $M$ be an irreducible orientable 3-manifold. Let $K$ be a knot in $M$.
Let $\mathcal{B}$ be a branched surface in $M$ disjoint from $K$. (see \cite{FO} \cite{O}
for definitions and facts about branched surfaces). Denote by
$N$ a fibered regular neighborhood of $\mathcal{B}$, by $\partial_h N$ the
horizontal boundary of $N$, and by $\partial_v N$ the vertical boundary of
$N$, as usual.

We say that a branched surface $\mathcal{B}$ is incompressible in $M-K$ if it
satisfies:

\begin{enumerate}
\item $\mathcal{B}$ has no discs of contact or half discs of contact.

\item $\partial_h N$ is incompressible and $\partial$-incompressible in
$(M-K)-int N$.

\item There are no monogons in $(M-K) - int N$.

We further say that $\mathcal{B}$ is meridionally
incompressible if:

\item $\partial_h N$ is meridionally incompressible, that is , there is
no disc $D$ in $M$, with $D\cap N =\partial D \subset \partial_h N$, so
that $D$ intersects $K$ transversely in one point.

We further say that $K$ is not parallel to $\mathcal{B}$ if:

\item $K$ is not parallel to $\partial_h N$, that is, there is no an annulus
$A$ in $M$, with $\partial A=A_0 \cup A_1$, so that $A_0=K$, and
$A\cap N = A_1 \subset \partial_h N$

\end{enumerate}

\begin{theorem}
Let $M$, $\mathcal{B}$, $K$ as above, with $\mathcal{B}$
incompressible.
\begin{enumerate}
\item Suppose $\mathcal{B}$
is meridionally incompressible. Then a surface carried with positive
weights by $\mathcal{B}$ is meridionally incompressible.
\item If $K$ is not parallel to $\mathcal{B}$, then $K$ is not parallel
to any  surface carried with positive
weights by $\mathcal{B}$.
\end{enumerate}
Then if $\mathcal{B}$ is meridionally incompressible and $K$ is not parallel
to it, any surface carried by $\mathcal{B}$ with positive weights is
quasi-Fuchsian.
\end{theorem}

\begin{proof} It is essentially the same proof as in Theorem 2 in \cite{FO},
with the obvious modifications.
\end{proof}

\bigskip
\centerline{\includegraphics{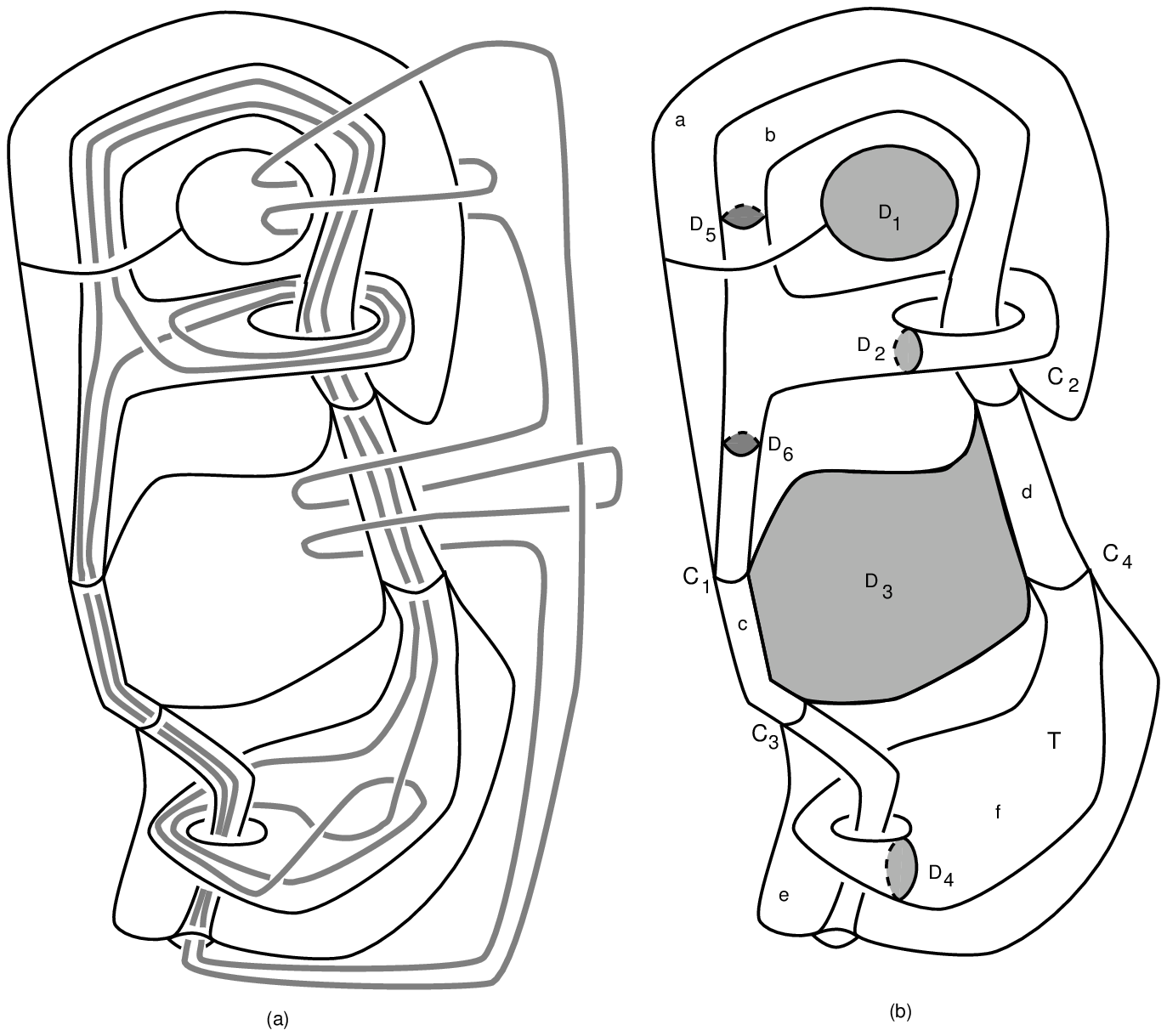}}
\centerline{\sr Figure 14}
\bigskip

Consider the knot $K$ and branched surface $\mathcal{B}$ shown in Figure 14(a). Note
that $\mathcal{B}$ has 4 singular curves, denoted $C_1,\ C_2,\ C_3,\ C_4$, as in
Figure 14(b). Note that $K$ is a 3-bridge knot. The knot $K$ is just
one in a collection of knots, to get more just make the knot to
intersect several times the discs $D_1,D_2,D_3,D_4$ shown in Figure 14(b).
But suppose that
$K$ intersects transversely the discs
$D_1,D_2,D_3,D_4$ in at least 2 points, that is, the minimal intersection
number of the knots with the discs, when isotoping the knot in the
complement of $\mathcal{B}$ is 2. Note that $K$ intersects the discs $D_5$ and
$D_6$ in exactly 2 points, because it is a 3-bridge knot. Suppose also
that the arc of the knot lying in the solid torus $T$ (shown in Figure
14(b)), is not parallel to $\partial T$; it is possible to do that, an
explicit example is in Figure 14(a).

The nonsingular part of $\mathcal{B}$ has six components, whose weights
$(a,b,c,d,e,f)$ are shown in Figure 14(b). Note that if we give the
weights $(1,2n-1,2n,2n-2,n,n-2)$, for $n\geq 3$, then this is a
collection of positive weights, which is consistent, and determines a
connected surface of genus $3n$.

If a knot $K$ is not hyperbolic then it is either a torus knot or a satellite
knot. Remember that by the classical work of Schubert, a satellite
3-bridge knot must be the connected sum of 2 two-bridge knots. It is known
that two-bridge knots do not contain any essential closed surface \cite{HT},
and from this it follows that the only essential surfaces in the connected
sum of 2 two-bridge knots are the swallow-follow tori. Also, torus knots do
not contain closed essential surfaces. This implies that a 3-bridge knot
which contains an essential surface of genus greater than 1 must be
hyperbolic.

\begin{theorem}
The surface $\mathcal{B}$ is meridionally
incompressible and $K$ is not parallel to it. So $K$ is a hyperbolic
3-bridge knot which contains quasi-Fuchsian surfaces of arbitrarily high
genus.
\end{theorem}

\begin{proof}[Sketch of proof] Let $N$ be a fibered neighborhood of $\mathcal{B}$. Note
that $S^3-N$ has 3 components, denoted by $N_1, N_2, N_3$, where say $N_3$
is the region that contains the knot, $N_1$ is the upper region, and $N_2$
the lower region.

Suppose that the part of $\partial_h N$ contained in $N_3$ is
compressible or meridionally compressible, and let $E$ be a compression or
meridian compression disc. Look at the intersections between $E$ and the
discs $D_1, D_3, D_5, D_6$. Let $\gamma$ be a simple closed curve of intersection which is innermost on $E$, so $\gamma$ bounds a disc
$E'\subset E$; suppose first that $E'$ is disjoint from $K$. The curve $\gamma$ also bounds a disc $D'$ in some $D_i$. Suppose $D'$ intersects $K$. If $D'$ is part of $D_1$ or $D_3$, then $K$ intersects the sphere $E'\cup D'$ several times always in the same direction, which is impossible. If $D'$ is part of $D_5$ or $D_6$ then it must intersect $K$ in two points, and then there is an arc of $K$ contained in the 3-ball bounded by $E'\cup D'$. But this implies that $K$ can be made disjoint from $D_2$ or $D_4$, or from $D_3$ or $D_1$, which is impossible by hypothesis. So $D'$ must be disjoint from $K$, and then an isotopy reduces the number of intersections between $E$ and the $D_i$. If $E'$ intersects $K$ once, then by a similar argument, $D'$ intersects $K$ also in a point, and then by an isotopy, we get a new compression disc with fewer intersections with the $D_i$. Suppose then that the intersection between $E$ and the $D_i$ consists only of arcs. Let $\gamma$ be an arc of intersection which is outermost on $E$, and which bounds a disc $E'$ disjoint from $K$. The arc $\gamma$ also bounds a disc $D'$ on some $D_i$. If $K$ is disjoint from $D'$, then by cutting $E$ with an outermost disc lying on $D'$ we get a new compression disc with fewer intersections with the $D_i$. If $K$ intersects $D'$ in one point, then it is not difficult to see that $K$ must intersect in one point one of $D_2$, $D_4$ or $D_1$, which is a contradiction. So if there is such a disc $E$, it must be disjoint from the $D_i$,  and by inspection it is not
difficult to check that such disc does not exist. The part of
$\partial_v N$
contained in $N_3$ consists of one annulus, corresponding to the curve
$C_2$. Again an innermost disc/outermost arc argument shows that there is
no monogon.

The part of $\partial_h N$ contained in $N_1$ consists
of a twice punctured genus two surface; it is not difficult to check that
it is incompressible. The part of $\partial_v N$ contained in $N_1$
consists of an annulus, corresponding to the curve $C_1$; it is also not
difficult to check that there is no monogon. Similarly, the part of
$\partial_h N$ contained in $N_2$ consists of a three punctured sphere and an once
punctured torus, and $\partial_v N$ consists of two annuli,
corresponding to the curves $C_3$ and $C_4$; again it is not
difficult to check that these are incompressible and that there is no
monogon.

To see that $K$ is not parallel to $\mathcal{B}$, suppose there is an annulus
$A$, with one boundary being $K$ and the other on $\mathcal{B}$. Again look at
the intersections between $A$ and the discs $D_i$, and get that the
arc of the knot that lies in the solid torus $T$ must be parallel to
$\partial T$, but this is not possible by the choice of such an  arc.
\end{proof}

The explicit knot shown in Figure 14(a) has more interesting properties,
it is a ribbon knot and it has unknotting number one, where a crossing
change is located in the arc contained in the solid torus $T$.

\begin{corollary}
There exist hyperbolic genus 3 closed
3-manifolds, in fact homology spheres, which contain incompressible
surfaces of arbitrarily high genus, so contain infinitely many
incompressible surfaces.
\end{corollary}

\begin{proof}
Let $K$ be a knot as in Theorem 6. Let $K(r)$ be the
manifold obtained by performing Dehn surgery on $K$ with slope $r$. If
$\Delta (r,\mu) > 1$, where $\mu$ denotes a meridian of $K$, then $K(r)$
is irreducible by \cite{GL1}, and $\mathcal{B}$ remains incompressible in $K(r)$ by
\cite{W}, for $K$ is not parallel to $\mathcal{B}$. If $\Delta (r,\mu) >2$, then
$K(r)$ is atoroidal by \cite{GL2}. So if $\Delta (r,\mu) >2$, $K(r)$ is an
atoroidal Haken manifold, hence it is hyperbolic. $K$ is a tunnel number
2 knot, hence each $K(r)$ has Heegaard genus $\leq 3$. Finally note that
among the $K(r)$ many are homology spheres.
\end{proof}

\begin{corollary}
There exist genus 2 closed 3-manifolds
which contain incompressible surfaces of
arbitrarily high genus, so they contain infinitely many incompressible
surfaces.
\end{corollary}

\begin{proof}
Let $K$ be a knot as in Theorem 6. Let $\Sigma (K)$ denote
the double cover of $S^3$ branched along $K$. As $K$ is a 3-bridge knot,
$\Sigma(K)$ has Heegaard genus 2. If $S$ is a surface carried by
$\mathcal{B}$ with positive weights, then as it is meridionally incompressible,
it lifts in $\Sigma(K)$ to a (possible disconnected) incompressible
surface \cite{GLi}.
\end{proof}

\begin{remark} It should be possible to say that the manifolds obtained in this corollary are hyperbolic; this will be the case if it is shown that the knots $K$ do not admit a tangle decomposing sphere.
\end{remark}

\end{document}